\documentclass[a4paper]{article}
\usepackage[dvips]{graphicx}
\usepackage{amssymb}
\usepackage{amsmath}
\usepackage{caption}
\usepackage{hyperref}

\makeatletter
\newcommand*{\centerfloat}{%
  \parindent \z@
  \leftskip \z@ \@plus 1fil \@minus \textwidth
  \rightskip\leftskip
  \parfillskip \z@skip}
\makeatother

\newcommand{\bs}{\boldsymbol}
\newcommand {\aplt} {\ {\raise-.5ex\hbox{$\buildrel<\over\sim$}}\ }
\makeatletter
\def\dddots{\mathinner{\mkern1mu\raise\p@
    \hbox{.}\mkern2mu\raise4\p@\hbox{.}\mkern2mu
    \raise7\p@\vbox{\kern7\p@\hbox{.}}\mkern1mu}}%
\makeatother
\newtheorem{theorem}{Theorem}
\newtheorem{Proposition}{Proposition}

\newtheorem{Corollary}{Corollary}

\newenvironment{AMS}{\small\bf 2010 AMS subject classification: }{}

\newfont{\BBB}{msbm10 scaled \magstep1}
\newfont{\BBS}{msbm10}

\begin{document} 

\title{\bf Optimal polynomial meshes and 
Caratheodory-Tchakaloff submeshes\\on the sphere\thanks{Work partially 
supported by the
DOR funds and the biennial project
CPDA143275 of the University of
Padova, and by the GNCS-INdAM.}}

\author{Paul Leopardi$^1$, Alvise Sommariva$^2$ and Marco Vianello$^2$}

\maketitle

\footnotetext[1]{Australian Government - Bureau of Meteorology, 
Docklands, 
Victoria, Australia\\e-mail: paul.leopardi@gmail.com}

\footnotetext[2]{Department of Mathematics, University of
Padova, Italy\\e-mail: alvise,marcov@math.unipd.it}

\begin{abstract} 
Using the notion of Dubiner distance, we give an elementary 
proof of the fact that good covering point configurations on the 2-sphere   
are optimal polynomial meshes. From these we extract 
Caratheodory-Tchakaloff (CATCH) submeshes for 
compressed Least Squares 
fitting.    
\end{abstract}
\vskip0.3cm
\noindent
\begin{AMS}
{\rm 41A10, 65D32.}
\end{AMS}

\section{Dubiner distance and polynomial meshes} 
In this note we focus on two notions that have played a relevant role 
in the theory of multivariate polynomial approximation during the last 20 
years: the notion of {\em polynomial mesh\/} and the notion of  
{\em Dubiner distance\/} in a compact set or manifold 
$K\subset \mathbb{R}^d$. Moreover, we connect the theory of polynomial 
meshes 
with the recent method of {\em Caratheodory-Tchakaloff 
(CATCH) subsampling\/}, working in 
particular on the sphere $S^2$. 

In what follows we denote by $\mathbb{P}^d_n(K)$ the subspace of 
$d$-variate polynomials of total degree not exceeding $n$ restricted to 
$K$, and by $N=N_n(K)=\mbox{dim}(\mathbb{P}^d_n(K))$ its dimension. 
For example we have that $N={n+3 \choose 3}=(n+1)(n+2)(n+3)/6$ for the 
ball in $\mathbb{R}^3$ and $N=(n+1)^2$ for the sphere $S^2$.  

We briefly recall that a polynomial mesh on $K$ 
is a sequence of finite norming subsets $\mathcal{A}_n\subset K$ such that 
the following 
polynomial inequality holds
\begin{equation} \label{wam}
\|p\|_{L^\infty(K)}\leq 
C\,\|p\|_{\ell^\infty(\mathcal{A}_n)}\;,\;\forall
p\in \mathbb{P}_n^d(K)\;,
\end{equation}
where 
$M_n=\mbox{card}(\mathcal{A}_n)=\mathcal{O}(N^\beta)$, $\beta\geq 1$. 
Indeed, since $\mathcal{A}_n$ is automatically
$\mathbb{P}_n^d(K)$-determining (i.e., polynomials vanishing there vanish 
everywhere on $K$), then $M_n\geq
N=\mbox{dim}(\mathbb{P}_n^d(K))=\mbox{dim}(\mathbb{P}_n^d(\mathcal{A}_n))$. 
Such a mesh is termed {\em optimal\/} when $\beta=1$.

In the case where $C$ is substituted by a sequence 
$C_n$ that increases subexponentially, 
\begin{equation} \label{wam2}
\|p\|_{L^\infty(K)}\leq
C_n\,\|p\|_{\ell^\infty(\mathcal{A}_n)}\;,\;\forall
p\in \mathbb{P}_n^d(K)\;,
\end{equation}
in particular when 
$C_n=\mathcal{O}(n^s)$, $s\geq 0$, we speak of a {\em
weakly admissible
polynomial mesh\/}. All these notions can be given 
for 
$K\subset \mathbb{C}^d$ but we restrict here to the real case. 
The notion of polynomial mesh was introduced in the seminal paper 
\cite{CL08} and then used from both the theoretical 
and the computational point of view, cf. e.g. 
\cite{BBCL12,BCLSV11,BDMSV11,BV12,K11,P12} 
and the references therein.

Polynomial meshes have indeed interesting computational features, e.g.
they 
\begin{itemize}

\item are affinely invariant and are stable under small 
perturbations \cite{PV13}

\item can be extended by algebraic transforms, finite
union and product \cite{BDMSV11,CL08}

\item contain computable near optimal
interpolation sets \cite{BCLSV11,BDMSV10}

\item are near optimal for uniform Least Squares (LS) approximation 
(cf. \cite[Thm. 1]{CL08}), 
namely
\begin{equation} \label{LSnorm}
\Lambda(\mathcal{A}_n)=\|\mathcal{L}_{\mathcal{A}_n}\|
=\sup_{f\in C(K),f\neq 0}\frac{\|\mathcal{L}_{\mathcal{A}_n}^wf\|_{L^\infty(K)}}
{\|f\|_{L^\infty(K)}}\leq C\,\sqrt{M_n}\;, 
\end{equation}
where $\mathcal{L}_{\mathcal{A}_n}$ is the $\ell^2(\mathcal{A}_n)$-orthogonal
projection operator
$C(K)\to \mathbb{P}_n^d(K)$ (the discrete LS operator at $\mathcal{A}_n$), 
from which easily follows    
\begin{equation} \label{LSerr}
\|f-\mathcal{L}_{\mathcal{A}_n}f\|_{L^\infty(K)}\leq \left(1+C\,\sqrt{M_n}\right)\,
\min_{p\in \mathbb{P}_{n}^d(K)}\|f-p\|_{L^\infty(K)}\;.
\end{equation}
\end{itemize}

We turn now to the notion of Dubiner distance in a compact set or 
manifold. Such a distance is defined as
\begin{equation}
\mbox{dist}_D(x,y)=\sup_{deg(p)\geq 1}\left\{\frac{1}{deg(p)}\,
|\cos^{-1}(p(x))-\cos^{-1}(p(y))|
\right\}\;,
\end{equation}
where the sup is taken over the polynomials $p\in 
\mathbb{P}_n^d(K)$ such that $\|p\|_{L^\infty(K)}\leq 1$. 

Introduced by M. Dubiner in the seminal paper \cite{D95}, it 
belongs to a family of three distances (the 
other two are the Markov 
distance and the Baran distance) that play an important role in 
multivariate 
polynomial approximation and have deep connections with multivariate 
polynomial inequalities. We refer the readers, e.g., to 
\cite{BLW04,BLW08,BLW08-2} 
and to the references therein for relevant properties and results.

As far as we know, 
until now the Dubiner distance 
is known analytically only in very few instances: the interval (where it 
coincides with the usual distance $|\cos^{-1}(x)-\cos^{-1}(y)|$), the 
cube, 
the simplex, the ball, and the sphere. All these cases are treated 
in \cite{BLW04}. In particular, it can be proved via the classical 
van der Corput-Schaake inequality that on the sphere it coincides with 
the usual {\em geodesic distance\/}, namely 
\begin{equation} \label{geodesic}
\mbox{dist}_D(x,y)=\gamma(x,y)=\cos^{-1}(\langle x,y\rangle)\;,\;\;\forall 
x,y\in S^2\;, 
\end{equation}
where $\langle x,y\rangle$ denotes the Euclidean scalar product in 
$\mathbb{R}^3$. 
 
A simple connection of the Dubiner distance with the theory of polynomial 
meshes is given by the following: 
\begin{Proposition}
Let $\mathcal{A}_n$ be a compact subset of a compact set or manifold 
$K\subset 
\mathbb{R}^d$  
whose covering radius with respect to the Dubiner distance does not exceed 
$\theta/n$, where $\theta\in (0,1)$ and $n\geq 1$, i.e.
\begin{equation} \label{covering}
\forall x\in K\;\exists y\in 
\mathcal{A}_n:\;\mbox{dist}_D(x,y)\leq \frac{\theta}{n}\;.
\end{equation}

Then, the following inequality holds
\begin{equation} \label{dubam}
\|p\|_{L^\infty(K)}\leq 
\frac{1}{1-\theta}\,\|p\|_{L^\infty(\mathcal{A}_n)}\;,\;\;\forall p\in 
\mathbb{P}_n^d(K)\;.
\end{equation}
\end{Proposition}
\vskip0.3cm
\noindent
In view of (\ref{covering}), the proof of Proposition 1 is an immediate 
consequence of the 
elementary inequality
\begin{equation} \label{dubinequal}
|p(x)|\leq |p(y)|+|p(x)-p(y)|\leq 
|p(y)|+n\,\mbox{dist}_D(x,y)\,\|p\|_{L^\infty(K)}\;.
\end{equation}
Observe, in particular, that $\mathcal{A}_n$ need not be discrete. 
In the case where (\ref{covering}) is satisfied by a sequence of finite 
subsets with $\mbox{card}(\mathcal{A}_n)=\mathcal{O}(N^\beta)$, 
$\beta\geq 1$, then 
these subsets clearly form a polynomial mesh like 
(\ref{wam}) for $K$, with $C=1/(1-\theta)$. 

Let us now focus on the case of the sphere, $K=S^2$. We recall that 
a sequence of finite point configurations $X_M\subset S^2$, with 
cardinality $M\geq 2$, is termed 
a ``good covering'' 
of the sphere if its {\em covering radius\/} 
\begin{equation} \label{covradius}
\eta(X_M)=\max_{x\in S^2}\min_{y\in X_M}|x-y|
\end{equation}
satisfies the inequality
\begin{equation} \label{estradius}
\eta(X_M)\leq \frac{\alpha}{\sqrt{M}}\;,
\end{equation}
for some $\alpha>0$ (cf., e.g., the survey paper \cite{HMS16}). It is then 
easy to prove the 
following result (that can be obtained also via a tangential Markov 
inequality on the sphere with exponent 1, cf. \cite{JSW98,JSW99}):
\begin{Proposition}
Let $\{X_M\}$, $M\geq 1$, be a good covering of $S^2$. Then 
for every fixed $\theta\in (0,1)$ the sequence $\mathcal{A}_n=X_{M_n}$, 
with 
\begin{equation} \label{Mn}
M_n=\lceil \sigma_n^2n^2\rceil\;,\;\;
\sigma_n=\frac{2\pi\alpha}{\theta(2\pi-\theta/n)}\sim 
\frac{\alpha}{\theta}\;,\;\;n\to \infty\;.
\end{equation}
is an optimal 
polynomial mesh of $S^2$ with $C=1/(1-\theta)$. 
\end{Proposition}
\vskip0.3cm
\noindent
{\em Proof}. By (\ref{covradius})-(\ref{estradius}) and simple geometric 
considerations, for 
every $x\in S^2$ there exists $y\in X_M$ such that we have the 
estimate
$$
\gamma(x,y)=2\sin^{-1}\left(\frac{|x-y|}{2}\right)
\leq 2\sin^{-1}\left(\frac{\alpha}{2\sqrt{M}}\right)
$$
(provided that $\alpha/(2\sqrt{M})\leq 1$ that is
$\sqrt{M}\geq \alpha/2$), 
where $\gamma$ is the geodesic distance, i.e., the Dubiner distance.   
By Proposition 1, in order to determine $M_n$ it is the sufficient to 
fulfill the 
inequality 
$$
2\sin^{-1}\left(\frac{\alpha}{2\sqrt{M}}\right)\leq \frac{\theta}{n}
$$
or equivalently 
\begin{equation} \label{intermediate}
\frac{\alpha}{2\sqrt{M}}\leq \sin\left(\frac{\theta}{2n}\right)\;.
\end{equation}
By the trigonometric inequality $\sin(t)\geq t(1-t/\pi)$, valid for $0\leq 
t\leq \pi$, 
we get
$$
\sin\left(\frac{\theta}{2n}\right)\geq 
\frac{\theta}{2n}\,\left(1-\frac{\theta}{2\pi n}\right)\;,
$$
and thus (\ref{intermediate}) is satisfied if 
$$
\frac{\alpha}{2\sqrt{M}}\leq\frac{\theta}{2n}\,
\left(1-\frac{\theta}{2\pi n}\right)\;,
$$
i.e. for $M\geq \sigma_n^2n^2$. $\;\;\;\;\;\square$
\vskip0.3cm
 
Among good covering configurations of $S^2$, an important role is played 
by the so-called ``quasi-uniform'' ones, that are those configurations 
with bounded ratio between the covering radius and the point separation   
(mesh ratio). Indeed, quasi-uniform configurations provide discretizations 
of the sphere that keep a low information redundancy. In \cite{HMS16} 
several quasi-uniform configurations are listed and their properties 
discussed. An interesting instance is the {\em zonal equal area 
configurations\/} (generated by zonal equal area partitions of the 
sphere), that turn out to be both, quasi-uniform and 
equidistributed in the sense of the surface area measure. In particular, 
they are theoretically good covering with $\alpha=3.5$ (but the 
numerical experiments suggest $\alpha=2.5$, cf. \cite{L06,L07}),  
and can be efficiently computed by the Matlab toolbox \cite{L06-2}.  
For example, taking $\theta=1/2$ by Proposition 2 we 
have the following quantitative result:  

\begin{Corollary} The zonal equal area configurations with 
$M_n=\left\lceil 49 \left(n-\frac{1}{4\pi}\right)^2 \right\rceil$ points 
are an 
optimal polynomial 
mesh of the sphere with $C=2$.  
\end{Corollary}
\vskip0.3cm
\noindent
The situation above is typical: polynomial meshes, even optimal ones, have 
often a large cardinality (though optimal as order of growth in $n$), 
which means large samples in the  
applications, for example in polynomial Least Squares. To this respect, 
it is useful to seek weakly admissible meshes with lower cardinality. 
This is exactly what we are going to do in the next Section, by the 
method of Caratheodory-Tchakaloff subsampling.  

To conclude this Section, we observe that Proposition 1 can be used to 
generate optimal polynomial meshes in 
other cases where the Dubiner distance is explicitly known. For example, 
in the square $K=[-1,1]^2$ 
\begin{equation} \label{dubsquare}
\mbox{dist}_D(x,y)=\max\left\{|\cos^{-1}(x_1)-\cos^{-1}(x_2)|,
|\cos^{-1}(y_1)-\cos^{-1}(y_2)|\right\}\;.
\end{equation}

We can then construct optimal polynomial meshes in $[-1,1]^2$ by the 
so-called {\em Padua points\/}. 
We recall that the Padua points of degree $k$ for 
the square $[-1,1]^2$ are given by the union of two Chebyshev-Lobatto
subgrids
\begin{equation} \label{padua}
\mathcal{P}_k=\mathcal{C}_{k+1}^{even}\times
\mathcal{C}_{k+2}^{odd} \cup \mathcal{C}_{k+1}^{odd}\times
\mathcal{C}_{k+2}^{even}\;,
\end{equation}
where
$\mathcal{C}_{s+1}=\{\cos(j\pi/s)\,,\,\,0\leq j\leq s\}$, and the
supscripts mean that only even or odd indexes are considered.
They are a near-optimal point set for polynomial interpolation of 
total degree $k$,
with a Lebesgue constant increasing as log square of the degree; 
cf., e.g., \cite{BCDMVX06}. Extending a similar result for univariate 
Chebyshev-like points, we can prove the following 
\begin{Proposition}
The Padua points $\mathcal{P}_{\nu_n}$,  
$\nu_n=\left\lceil\frac{\pi n}{\theta}\right\rceil$, $0<\theta<1$ 
(cf. (\ref{padua})),  
are an
optimal polynomial
mesh of the square with $C=(1-\theta)^{-1}$ and cardinality 
$M_n=(\nu_n+1)(\nu_n+2)/2$.
\end{Proposition}
\vskip0.3cm
\noindent
{\em Proof}. 
In view of (\ref{dubsquare}), by simple geometric considerations we 
get easily that 
$\mbox{dist}_D(x,\mathcal{P}_{\nu_n})
\leq \pi/\nu_n\leq \theta/n$, for every $x\in [-1,1]^2$. The result then 
follows immediately by Proposition 1. $\;\;\;\;\;\square$
\vskip0.3cm 

We stress that the cardinality is asymptotically $\nu_n^2/2$, that is 
essentially half the cardinality obtainable by embedding the problem 
in the tensorial polynomial space $\mathbb{P}_n^1\otimes \mathbb{P}_n^1$,  
and using for example product Chebyshev points of degree $\nu_n$ (cf. 
\cite{BV12,V14}).

\subsection{Caratheodory-Tchakaloff (CATCH) submeshes}
In order to reduce the cardinality of a polynomial mesh, a feature that is 
relevant in applications, we may try to relax the boundedness 
requirement for 
the ratio $\|p\|_{L^\infty(K)}/\|p\|_{\ell^\infty(\mathcal{A}_n)}$, 
seeking a {\em weakly admissible mesh\/} contained in the original 
one, 
where the ratio is allowed to increase algebraically with 
respect to 
the polynomial space dimension $N$. 

In principle, this can be done by 
computing Fekete points $\mathcal{F}_n\subset \mathcal{A}_n$. 
These 
are a subset of $\mathcal{A}_n$ that maximizes the Vandermonde 
determinant (not unique in general), so that the 
corresponding cardinal polynomials are bounded by 1 in 
$\ell^\infty(\mathcal{A}_n)$, and by $C$ in $L^\infty(K)$. This entails 
that such discrete Fekete points are unisolvent and form a weakly 
admissible mesh of cardinality $N$, $C_n=CN$ being a bound on their 
Lebesgue constant \cite{CL08}. Unfortunately, even the discrete Fekete 
points (as 
the continuous one) are difficult and costly to compute. In several 
papers, approximate Fekete points extracted from polynomial meshes 
have been computed by greedy algorithms 
based on standard numerical linear algebra routines; cf., e.g., 
\cite{BCLSV11,BDMSV10,BDMSV11}. 
These points 
work effectively for interpolation, and are asymptotically distributed 
as the continuous Fekete points, but no rigorous bound has been 
proved 
for their Lebesgue constant. 

In the case of the sphere, the continuous 
Fekete (maximum determinant) points have been computed by a difficult 
numerical nonconvex 
optimization 
up to degrees in the hundreds, estimating also numerically the 
corresponding Lebesgue 
constants; cf. \cite{WS01}. As it is well-known, the difficulties 
of polynomial interpolation on the sphere have led to the alternative 
approach of {\em hyperinterpolation\/}, cf. the seminal paper \cite{S95} 
and the numerous developments in the following 20 years. 

In the present note, starting from optimal polynomial meshes 
of the sphere, we explore an alternative discrete approach, 
that can be considered a sort of fully discrete hyperinterpolation, 
namely the extraction of {\em Caratheodory-Tchakaloff submeshes\/}.
These are computable by Linear or Quadratic Programming, and 
there are rigorous bounds for the corresponding constants $C_n$. 

First, we recall a discrete version of the Tchakaloff theorem, a 
cornerstone 
of quadrature theory, whose proof is based on 
the Caratheodory theorem
about finite dimensional conic combinations (cf., e.g., \cite{BV10}). 
We focus here on 
total-degree polynomial 
spaces and we recall the proof to exhibit the connection with 
the Caratheodory 
theorem.  
\begin{theorem}
Let $\mu$ be a multivariate discrete measure
supported at a finite set
$X=\{x_i\}\subset \mathbb{R}^d$,
with correspondent
positive weights (masses)
$\bs{\lambda}=\{\lambda_i\}$, $i=1,\dots,M$.

Then, there exists a quadrature formula with nodes $T=\{t_j\}\subseteq X$
and positive
weights $\bs{w}=\{w_j\}$, $1\leq j\leq m\leq 
N_\nu=\mbox{{\em dim\/}}(\mathbb{P}_\nu^d(X))$, such that
\begin{equation} \label{qformula}
\int_X{p(x)\,d\mu}=\sum_{i=1}^M{\lambda_i\,p(x_i)}
=\sum_{j=1}^m{w_j\,p(t_j)}\;,\;\;\forall p\in \mathbb{P}_\nu^d(X)\;.
\end{equation}
\end{theorem}
\vskip0.3cm
\noindent
{\em Proof\/} (cf., e.g., \cite{PSV16}). Let $\{p_1,\dots,p_{N_\nu}\}$ be 
a 
basis 
of $\mathbb{P}_\nu^d(X)$, 
and
$V=(v_{ij})=(p_j(x_i))$ the Vandermonde-like matrix of the basis 
computed at the support points. If $M>N_\nu$ (otherwise there is
nothing to prove), existence of a positive quadrature formula
for $\mu$ with cardinality not exceeding $N_\nu$ can be immediately 
translated
into existence of a nonnegative solution with at most $N_\nu$ nonvanishing
components to the underdetermined linear system
\begin{equation} \label{momeqs}
V^t\bs{u}=\bs{b}\;,\;\;\bs{u}\geq \bs{0}\;,
\end{equation}
where
\begin{equation} \label{moments}
\bs{b}=V^t\bs{\lambda}=\left\{\int_X{p_j(x)\,d\mu}\right\}\;,\;\;1\leq
j\leq N_\nu\;,
\end{equation}
is the vector of $\mu$-moments of the basis $\{p_j\}$.

Existence then holds by the
well-known Caratheodory theorem applied to the columns of $V^t$, which  
asserts that a conic (i.e., with positive coefficients) combination
of any number of vectors in $\mathbb{R}^N$ can be rewritten as a conic 
combination of a linearly independent subset of at most $N_\nu$ of 
them.$\;\;\;\;\;\square$
\vskip0.3cm

We may term $T=\{t_j\}$ a set of 
{\em Caratheodory-Tchakaloff (CATCH) quadrature points\/}. 
We apply now the Tchakaloff theorem to the extraction 
of a weakly admissible submesh from a polynomial mesh.
\begin{Proposition} 
Let $\mathcal{A}_n\subset K$ be a polynomial mesh like (\ref{wam}) 
for $K$  
with 
cardinality 
$M_n>N_{2n}=\mbox{dim}(\mathbb{P}_{2n}^d(K))$, let $\mu$ be the discrete 
measure with unit weights supported at $\mathcal{A}_n$,   
and let $T_{2n}=\{t_j\}$ be the $m\leq N_{2n}$ Caratheodory-Tchakaloff 
quadrature points 
for exactness degree $\nu=2n$ 
extracted from $\mathcal{A}_n$, with corresponding weights 
$\bs{w}=\{w_j\}$, 
$1\leq j\leq m$. Moreover, let $\mathcal{L}_{\mathcal{T}_{2n}}^wf \in 
\mathbb{P}_n^d(K)$ be 
the weighted discrete 
least squares polynomial on $T_{2n}$ for $f\in C(K)$.  

Then $T_{2n}$ is a weakly-admissible polynomial mesh like (\ref{wam2}) 
for $K$ with 
\begin{equation} \label{Cn}
C_n=C\sqrt{M_n} 
\end{equation}
(that we may term a 
Caratheodory-Tchakaloff submesh), and the following estimate holds for 
the corresponding weighted least squares approximation  
\begin{equation} \label{WLSerr} 
\|f-\mathcal{L}_{T_{2n}}^wf\|_{L^\infty(K)}\leq \left(1+C_n\right)\,
\min_{p\in \mathbb{P}_{n}^d(K)}\|f-p\|_{L^\infty(K)}\;.
\end{equation}
\end{Proposition}    
\vskip0.3cm
\noindent
{\em Proof} (cf., e.g., \cite{V16}).  
Since the Caratheodory-Tchakaloff quadrature is exact
in $\mathbb{P}_{2n}^d(K)$, we get the basic $\ell^2$-identity
\begin{equation} \label{2-norms}
\|p\|^2_{\ell^2(\mathcal{A}_n)}=
\sum_{i=1}^{M_n}{p^2(x_i)}=
\sum_{j=1}^m{w_j\,p^2(t_j)}   
=\|p\|^2_{\ell^2_{\bs{w}}(T_{2n})}\;,\;\;\forall p\in
\mathbb{P}_{n}^d(K)\;.
\end{equation}
Then we can write
$$
\|p\|_{L^\infty(K)}\leq C\,\|p\|_{\ell^\infty(\mathcal{A}_n)}
\leq C\,\|p\|_{\ell^2(\mathcal{A}_n)}
=C\,\|p\|_{\ell^2_{\bs{w}}(T_{2n})}
$$
\begin{equation} \label{catchwam}
\leq 
C\,\sqrt{\sum_{j=1}^m{w_j}}\;\;\|p\|_{\ell^\infty(T_{2n})}
=C\,\sqrt{M_n}\,\|p\|_{\ell^\infty(T_{2n})}\;,
\end{equation}
i.e., $T_{2n}$ is a weakly-admissible polynomial mesh 
for $K$ with $C_n=C\sqrt{M_n}$.

Concerning the weighted Least Squares polynomial approximation 
on $T_{2n}$, we recall that it is defined by
\begin{equation} \label{WLS}
\|f-\mathcal{L}_{T_{2n}}^wf\|_{\ell^2_{\bs{w}}(T_{2n})}
=\min_{p\in \mathbb{P}_{n}^d(K)}\|f-p\|_{\ell^2_{\bs{w}}(T_{2n})}
\end{equation}
and that $\mathcal{L}_{T_{2n}}^wf$ is a $\ell^2_{\bs{w}}(T_{2n})$-orthogonal 
projection, i.e., $f-\mathcal{L}_{T_{2n}}^wf$ is 
$\ell^2_{\bs{w}}(T_{2n})$-orthogonal to 
$\mathbb{P}_{n}^d(K)$. 
Then,  
$$
\|\mathcal{L}_{T_{2n}}^wf\|_{L^\infty(K)}\leq 
C\,\|\mathcal{L}_{T_{2n}}^wf\|_{\ell^\infty(\mathcal{A}_n)}
\leq C\,\|\mathcal{L}_{T_{2n}}^wf\|_{\ell^2(\mathcal{A}_n)}
$$
$$
=C\,\|\mathcal{L}_{T_{2n}}^wf\|_{\ell^2_{\bs{w}}(T_{2n})}
\leq C\,\|f\|_{\ell^2_{\bs{w}}(T_{2n})}
\leq C\,\sqrt{M_n}\,\|f\|_{L^\infty(K)}\;,
$$ 
that is 
\begin{equation} \label{WLSnorm}
\Lambda_w(T_{2n})=\|\mathcal{L}_{T_{2n}}^w\|=\sup_{f\in C(K),f\neq 
0}\frac{\|\mathcal{L}_{T_{2n}}^wf\|_{L^\infty(K)}}
{\|f\|_{L^\infty(K)}}\leq C_n=C\,\sqrt{M_n}\;,
\end{equation} 
from which (\ref{WLSerr}) easily follows.$\;\;\;\;\;\square$
\vskip0.3cm

Observe that the error estimate (\ref{WLSerr}) for weighted 
discrete Least 
Squares on the Caratheodory-Tchakaloff submesh, turns out to coincide 
with the natural error estimate (\ref{LSerr}) for 
unweighted Least Squares on the original polynomial mesh. 
In some sense, Caratheodory-Tchakaloff weighted Least 
Squares on the 
submesh catches all the relevant information from the polynomial mesh, as 
far as 
polynomial approximation in $\mathbb{P}_{n}^d(K)$ is concerned.  
We recall that the best uniform approximation error in 
$\mathbb{P}_{n}^d(K)$ can be estimated by 
the regularity of $f$, on compact sets $K$ admitting a Jackson-like 
theorem; cf., e.g., \cite{P09}.

We can now apply the results above to optimal polynomial meshes 
on the sphere. Indeed, from Proposition 2 and 4 and Corollary 1 and 2 
we get immediately the following
\begin{Corollary}
Let $\mathcal{A}_n$ be a good covering optimal polynomial mesh 
as in Proposition 2, 
and let $T_{2n}$ be the extracted Caratheodory-Tchakaloff submesh 
(with corresponding weights). 

Then, $T_{2n}$ is a weakly admissible mesh for the sphere with 
cardinality $N_{2n}=\mbox{{\em dim\/}}(\mathbb{P}^3_{2n}(S^2))=(2n+1)^2$, 
and (\ref{WLSerr}) holds for the corresponding weighted Least Squares 
polynomial approximation $\mathcal{L}_{T_{2n}}^wf$ to $f\in C(S^2)$, where 
\begin{equation} \label{Cnsphere}
C_n=\frac{\sigma_n\,n}{1-\theta}\sim \frac{\alpha n}{\theta (1-\theta)}
\;,\;\;n\to \infty\;.
\end{equation}
In particular, for a Caratheodory-Tchakaloff submesh of the zonal equal 
area configurations of Corollary 1, we have
\begin{equation} \label{CnEA}
C_n=\frac{14 n}{1-(4\pi n)^{-1}}\sim 14 n\;,\;\;n\to \infty\;.
\end{equation}

\end{Corollary}
\vskip0.3cm
\noindent
We observe that by ((\ref{LSnorm}), (\ref{WLSnorm}) and (\ref{Cnsphere}) 
we get $\mathcal{O}(n)$ estimates for the least squares operator norms, 
whereas the best projection 
operators on $\mathbb{P}^3_{n}(S^2)$ have a $\mathcal{O}(n^{1/2})$ norm; 
cf., e.g., \cite{SW00}. On the other hand, (\ref{WLSnorm}) turns out to 
be an overestimate of the actual norm, as we shall see in the numerical 
examples. 

\subsection{Computational issues and numerical examples}
In order to compute a sparse nonnegative solution to the underdermined 
system (\ref{momeqs})-(\ref{moments}), that exists by Tchakaloff theorem 
(Theorem 1 
above), there are a number of different approaches available. We focus 
here on the case of 
the sphere, where we use the classical spherical harmonics basis to 
define the Vandermonde-like matrix $V$.

A first approach resorts to Quadratic 
Programming, namely to the  
NonNegative Least Squares problem
\begin{equation} \label{NNLS}
\mbox{QP}:\;\;\;\left\{
\begin{array}{ll}
\min{\|V^t\bs{u}-\bs{b}\|_2}\\
\bs{u}\geq \bs{0}
\end{array}   
\right.
\end{equation}
which can be solved by the Lawson-Hanson active set method, 
which naturally seeks a sparse solution; cf. \cite{PSV16,SV15} and the 
references therein. The nonzero components of $\bs{u}$ identify the 
weights $\bs{w}=\{w_j\}$ and the corresponding CATCH submesh $T_{2n}$. 

A second approach is based on Linear Programming (cf. 
\cite{PSV16,RB15,T15}), namely
\begin{equation} \label{LP}
\mbox{LP}:\;\;\;\left\{\begin{array}{ll}
\min{\bs{c}^t\bs{u}}\\
V^t\bs{u}=\bs{b}\;,\;\;\bs{u}\geq
\bs{0}
\end{array}
\right.
\end{equation}
where the constraints identify a
polytope (the feasible region), and the vector $\bs{c}$ is suitably 
chosen (cf. \cite{PSV16,RB15}). Solving the problem by the classical 
Simplex Method, we get a vertex of the polytope, that is  
a nonnegative sparse solution to the underdetermined system. 

A third combinatorial approach (Recursive
Halving Forest), based on the SVD, is proposed in \cite{T15} 
and essentially applied to the reduction
of Cartesian tensor cubature measures. 

It is worth observing that sparsity cannot be ensured by the
standard Compressive Sensing approach to underdetermined systems, such as
the Basis Pursuit algorithm that minimizes $\|u\|_1$ (cf., e.g., 
\cite{FR13}), since 
$\|u\|_1=\sqrt{M_n}$ is here constant by construction (being the 
quadrature formula applied to the constant polynomial $p\equiv 1$).

In our Matlab codes for Caratheodory-Tchakaloff Least Squares we have 
adopted both the QP approach (via an optimized version of 
the {\tt lsqnonneg\/} function), and the LP approach (via the Simplex 
Method in the Matlab interface of the CPLEX package); cf. \cite{PSV16}.   

In Table 1, we report the numerical results corresponding to the 
extraction of CATCH submeshes from zonal equal area meshes of $S^2$, 
for a sequence of degrees. 
All the quantities are rounded to the first decimal digit. 
We have that 
the cardinality of the CATCH submeshes is   
$\mbox{dim}(\mathbb{P}_{2n}^3(S^2))=(2n+1)^2$, and that the 
Compression Ratio, $C_{ratio}=card(\mathcal{A}_n)/card(T_{2n})$,  
increases, approaching the asymptotic value $49/4=12.25$, cf. 
Corollary 1. Since $\sum w_j=card(\mathcal{A}_n)$, the average 
CATCH weight 
turns out to coincide with the Compression Ratio. Notice that the 
minimum of the CATCH weights computed by QP is much smaller than the 
minimum of the CATCH weights computed by LP. 

Moreover, we 
see that the compressed least squares operator norms 
$\Lambda_w(T_{2n})$ are close to 
the norm $\Lambda(\mathcal{A}_n)$ of the least squares operator 
on the starting mesh, 
with a slightly better behavior of CATCH submeshes extracted by QP 
with respect to those extracted by LP. On the other hand, all the norms  
are much lower than 
the theoretical overestimate $C_n\sim 14n$ in Corollary 2, having 
substantially a $\mathcal{O}(n^{1/2})$ increase (at least in the 
considered degree range). 

In Table 2, we report the reconstruction errors by Least Squares (in the 
$L^\infty(S^2)$-norm, numerically evaluated on a fine control grid), 
for three test functions with different degree of regularity 
$$
f_1(x_1,x_2,x_3)=(x_1+x_2+x_3)^{15}\;,\;\;
f_2(x_1,x_2,x_3)=\exp(x_1+x_2+x_3)/10\;,\;\;
$$
\begin{equation} \label{testf}
f_3(x_1,x_2,x_3)=(|x_1|+|x_2|+|x_3|)/10\;,\;\;
\end{equation}
namely a polynomial, a smooth function and a function with singular 
points on the sphere (the latter two taken from
\cite{WS01}). 

Finally, in Figure 1 we display the CATCH submesh
extracted by QP from a zonal equal area mesh for degree $n=5$.

\begin{table}
\begin{center}
\caption{Cardinality of zonal equal area meshes $\mathcal{A}_n$ and 
of their CATCH 
submeshes $T_{2n}$,  
Compression Ratio, weight ratios, discrete least squares operator norms 
(the CATCH weights and 
submeshes have been obtained by (\ref{NNLS}) via {\tt lsqnonneg\/} and 
by (\ref{LP}) via CPLEX).}
{\footnotesize
\begin{tabular}{c c c c c c c c}
deg $n$ & 2 & 5 & 8 & 11 & 14 & 17 & 20\\ 
\hline
$card(\mathcal{A}_n)$ & 181 & 1187 & 3074 & 5844 & 9496 & 14029 & 19445\\ 
\hline 
$card(T_{2n})$ & 25 & 121 & 289 & 529 & 841 & 1225 & 1681\\ 
\hline
$C_{ratio}$ ($=w_{avg}$) & $7.2$ & $9.8$ & $10.6$ & $11.0$ & $11.3$ & 
$11.5$ & $11.6$\\ 
\hline
QP: $w_{max}/w_{avg}$ & 2.2 & 2.6 & 2.5 & 2.5 & 2.6 & 2.4 & 
2.6\\
LP: $w_{max}/w_{avg}$ & 2.1 & 3.1 & 2.8 & 3.1 & 3.1 & 3.0 & 3.0\\
\hline
QP: $w_{min}/w_{avg}$ & 2.1e-2 & 5.4e-4 & 7.6e-5 & 8.8e-6 & 
2.4e-6 & 8.1e-6 & 2.6e-6\\
LP: $w_{min}/w_{avg}$ & 9.1e-2 & 4.3e-4 & 2.3e-3 & 6.2e-4 & 1.1e-3 
& 2.4e-3 & 1.3e-3\\
\hline
$\Lambda(\mathcal{A}_n)$ & 2.2 & 3.3 & 4.2 & 4.9 & 5.6 & 6.2 & 6.7\\ 
\hline
$1.5\,n^{1/2}$ & 2.1 & 3.4 & 4.2 & 5.0 & 5.6 & 6.2 & 6.7\\
\hline
QP: $\Lambda_w(T_{2n})$ & 2.5 & 3.7 & 4.6 & 5.4 & 6.0 & 6.5 & 7.1\\ 
\hline 
LP: $\Lambda_w(T_{2n})$ & 2.6 & 4.2 & 5.6 & 6.6 & 7.4 & 7.5 & 8.4\\ 
\hline
\end {tabular}
\/}
\end{center}
\end{table}

\begin{table}
\begin{center}
\caption{Absolute errors ($L^\infty(S^2)$-norm) in the reconstruction of 
the three 
test functions (\ref{testf}) by unweighted Least Squares on zonal equal 
area meshes $\mathcal{A}_n$ 
and by weighted Least Squares on their CATCH
submeshes $T_{2n}$.}
{\footnotesize
\begin{tabular}{c c c c c c c c c}
 & deg $n$ & 2 & 5 & 8 & 11 & 14 & 17 & 20\\
\hline
 & LS & 1.4e+5 & 4.1e+4 & 1.5e+4 & 5.3e+2 & 3.7e+1 & 
8.4e-10 & 9.1e-10\\
$f_1$ & CATCH$_{QP}$ & 1.7e+5 & 4.8e+4 & 1.4e+4 & 5.1e+2 & 3.7e+1 & 
8.4e-10 & 
6.4e-10\\
 & CATCH$_{LP}$ & 2.0e+5 & 4.6e+4 & 1.6e+4 & 6.5e+2 & 4.3e+1 & 5.8e-10 & 
6.7e-10\\
\hline
& LS & 1.3e-1 & 1.6e-3 & 5.1e-6 & 5.9e-9 & 3.3e-12 &
5.6e-15 & 6.2e-15\\
$f_2$ & CATCH$_{QP}$ & 1.7e-1 & 1.7e-3 & 5.1e-6 & 6.0e-9 & 3.3e-12 &
2.8e-15 & 1.9e-15\\ 
 & CATCH$_{LP}$ & 1.3e-1 & 1.7e-3 & 4.5e-6 & 6.7e-9 & 3.3e-12 &
2.1e-15 & 2.6e-15\\
\hline
& LS & 5.0e-1 & 3.2e-1 & 1.5e-1 & 1.4e-1 & 9.5e-2 &
9.2e-2 & 7.4e-2\\ 
$f_3$ & CATCH$_{QP}$ & 6.0e-1 & 3.4e-1 & 1.9e-1 & 1.6e-1 & 1.2e-1 & 
9.7e-2 & 7.4e-2\\  
 & CATCH$_{LP}$ & 6.0e-1 & 3.7e-1 & 1.9e-1 & 1.5e-1 & 1.3e-1 & 
1.0e-1 & 1.1e-1\\ 
\hline

\end {tabular}
\/}
\end{center}
\end{table}

\begin{figure}
\centerfloat
\includegraphics[scale=0.40,clip]{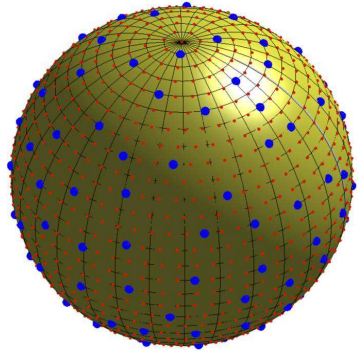}
\caption{CATCH submesh (121 points, bullets) extracted by NNLS from
a zonal equal area mesh
(1187 points, dots) of the sphere
for degree $n=5$.}
\end{figure}

\end{document}